\documentclass[12pt,reqno]{amsart}
\usepackage{amssymb}
\usepackage{amsfonts}
\usepackage{amssymb,latexsym}
\usepackage{enumerate}
\usepackage{mathrsfs}
\usepackage{url}
\usepackage{tikz}
\usetikzlibrary{calc, arrows.meta}

\makeatletter
\@namedef{subjclassname@2020}{%
  \textup{2020} Mathematics Subject Classification}
\makeatother

  [2002/01/19 v2.2g %
    AMS font definitions%
  ]
\DeclareFontFamily{U}{euf}{}
\DeclareFontShape{U}{euf}{m}{n}{%
  <5><6><7><8><9>gen*eufm%
  <10><10.95><12><14.4><17.28><20.74><24.88>eufm10%
  }{}
\DeclareFontShape{U}{euf}{b}{n}{%
  <5><6><7><8><9>gen*eufb%
  <10><10.95><12><14.4><17.28><20.74><24.88>eufb10%
  }{}

  [2002/01/19 v2.2g %
    AMS font definitions%
  ]
\DeclareFontFamily{U}{msb}{}
\DeclareFontShape{U}{msb}{m}{n}{%
  <5><6><7><8><9>gen*msbm%
  <10><10.95><12><14.4><17.28><20.74><24.88>msbm10%
  }{}

  [2002/01/19 v2.2g %
    AMS font definitions%
  ]
\DeclareFontFamily{U}{msa}{}
\DeclareFontShape{U}{msa}{m}{n}{%
  <5><6><7><8><9>gen*msam%
  <10><10.95><12><14.4><17.28><20.74><24.88>msam10%
  }{}

\newtheorem{theorem}{Theorem}[section]

\theoremstyle{definition}
\newtheorem{remark}[theorem]{Remark}

\numberwithin{equation}{section} 

\begin{document}

\title[Reciprocity law]
{On the reciprocity law  in $\mathbb{F}_{q}[t]$}
\author{Su Hu}
\address{Department of Mathematics, South China University of Technology, Guangzhou, Guangdong 510640, China}
\email{mahusu@scut.edu.cn}

\author{Enci Wang}
\address{Department of Mathematics, South China University of Technology, Guangzhou, Guangdong 510640, China}
\email{maecw@mail.scut.edu.cn}






\begin{abstract}
In 1991, Rousseau gave a new proof of Gauss's quadratic reciprocity by comparing two distinct coset representations of the group $(\mathbb{Z}_{p}^{*} \times \mathbb{Z}_{q}^{*}) / U$ using the Chinese Remainder Theorem, without Gauss's Lemma. In this paper, we extend Rousseau's approach to $\mathbb{F}_{q}[t]$, providing a new, elementary proof of the reciprocity law for the $d$th power residue symbol, where $d$ is any divisor of $q-1$. \end{abstract}
 \maketitle

\section{Introduction}
Let $\left(\frac{p}{q}\right)$ be the Legendre symbol.  Gauss's reciprocity law reads as follows
\begin{equation*}
\left(\frac{p}{q}\right)\left(\frac{q}{p}\right)=(-1)^{\frac{p-1}{2}\frac{q-1}{2}}
\end{equation*}
for $p, q$ are distinct odd primes. In the history, there are many methods to prove this.
One elegant proof  is due to Rousseau (1991) (see \cite{Rousseau}). Let $\mathbb{Z}_{m}^{*}$ denote the multiplicative group of reduced residues modulo $m$. Rousseau's approach is based on comparing two distinct coset representations of a group constructed via the Chinese remainder theorem
\[
G = (\mathbb{Z}_{p}^{*} \times \mathbb{Z}_{q}^{*}) / U \cong \mathbb{Z}_{pq}^{*} / U,
\]
where \( U = \{(1,1), (-1,-1) \} \). The first representation consists of the elements
\begin{equation*}
\{(i,j): i=1,2, \ldots, p-1; \; j=1,2,\ldots, (q-1)/2\},
\end{equation*}
while the second representation is given by
\begin{equation*}
\{(k \bmod p, k \bmod q): k=1,2, \ldots, (pq-1)/2, \gcd(k, pq)=1\}.
\end{equation*}
Rousseau's method reveals clearly that the reciprocity law is a direct consequence of the Chinese remainder theorem.

In this paper, let $q$ be a power of an odd prime number $p$, $A=\mathbb{F}_{q}[t]$ be the polynomial 
ring over $\mathbb{F}_{q}$ and $A_{+}$ be the set of monic polynomials in $A$.
Let $d$ be any divisor of $q-1$ and $P\in A$ be an irreducible polynomial. For $a\in A$ and $P$ dose not divide $a$, the
$d$th power residue symbol in $A$ is defined in analogy with the Euler criterion (see \cite[p. 24]{Rosen}):
\begin{equation}\label{def}
\left(\frac{a}{P}\right)_{d} \equiv a^{\frac{|P|-1}{d}} \pmod{P},
\end{equation}
where $|P|=q^{\textrm{deg}P}$ and denotes $\left(\frac{a}{P}\right)_{d}=0$ if $P\mid a$.
Let $P(t), Q(t)\in A$ be distinct  irreducible monic polynomials, there exists  a reciprocity law for the $d$th power residue symbol 
(see \cite[p. 25, Theorem 3.3]{Rosen}). 
\begin{theorem}
\begin{equation}\label{(1)}
\left(\frac{P}{Q}\right)_{d}\left(\frac{Q}{P}\right)_{d}^{-1}=(-1)^{\frac{q-1}{d}\textrm{deg}P\textrm{deg}Q}.
\end{equation}
\end{theorem}
\begin{remark} 
In contrast to the complex development of the higher reciprocity laws over number fields, the history for the polynomial ring $\mathbb{F}_{q}[t]$ is more direct.
Dedekind first proved an analogue of the quadratic reciprocity law in this setting.  The general $d$th reciprocity law for $\mathbb{F}_{q}[t]$ was established by  P. K. Schmidt in 1928.
Subsequently, Carlitz independently discovered the theorem and provided several elegant proofs, demonstrating that its formulation and proof do not require the extensive machinery such as the class field theory.
(See \cite[p. 23]{Rosen}).
\end{remark}
In this note, we shall generalize Rousseau's method to give a new proof of (\ref{(1)}).
Because the unit group of $\mathbb{Z}$ is $\{\pm 1\}$, whereas that of $A = \mathbb{F}_{q}[t]$ is the larger group $\mathbb{F}_{q}^{*}$, generalizing this method requires a more refined treatment.

 \section{Proof}
Let $g$ be a generator of $\mathbb{F}_{q}^{*}$, and denote $\eta=g^{\frac{q-1}{d}}$, an element of order $d$ in $\mathbb{F}_{q}^{*}$. Let
$$
G_{1}=\left(A/(P(t))\right)^{*}\times\left(A/(Q(t))\right)^{*},
$$
$U=\left\langle (\eta,\eta)\right\rangle$
be the cyclic group of $G_{1}$ generated by $(\eta,\eta)$ and $G=G_{1}/U$.
Let 
\begin{equation}\label{(Sp)}
S_{P}:=\{f\in A: 0\leq \textrm{deg}f <\textrm{deg}P\}=\{f\in A: 0<|f|<|P|\}.
\end{equation}
Let $S_{Q}$ be the coset representation of the group  $\left(A/(Q(t))\right)^{*}/\langle \eta \rangle$,
we have 
  \begin{equation}\label{(2)}
  S_{Q}\cup \eta S_{Q}\cup \cdots \cup\eta^{d-1}S_{Q}=\{h\in A: 0<|h|< |Q|\}
  \end{equation} 
 and $\# S_{Q}=\frac{|Q|-1}{d}$.
 Denote  $S_{1}=S_{P}\times S_{Q}$. It is seen that $S_{1}$ is a representation for the cosets of $U$ in $G$, and the product of all the
 elements in $G$ under the representation $S_{1}$ is 
\begin{equation}\label{(3)}
\pi = \left(\left( \prod_{f \in S_{P}} f \right)^{\frac{|Q|-1}{d}}, \left(\prod_{h \in S_{Q}} h\right)^{|P|-1}\right)U.
\end{equation}  
Then by (\ref{(2)}) we have
  $$\left(\prod_{h \in S_{Q}} h\right)^{d}\eta^{\#S_{Q}\frac{d(d-1)}{2}}=\left(\prod_{h \in S_{Q}} h\right)^{d}\eta^{\frac{|Q|-1}{d}\frac{d(d-1)}{2}}=\prod_{\substack{h \in A\\0<|h|< |Q|}}h.$$
Since  $$\eta^{\frac{d(d-1)}{2}}=g^{\frac{q-1}{d}\frac{d(d-1)}{2}}=(-1)^{d-1},$$
we get 
 $$\left(\prod_{h \in S_{Q}} h\right)^{d}=(-1)^{\frac{|Q|-1}{d}}\prod_{\substack{h \in A\\0<|h|< |Q|}}h$$
 by noticing that $|Q|=q^{\textrm{deg} Q}$ is odd. Thus 
 $$\left(\prod_{h \in S_{Q}} h\right)^{|P|-1}=\left(\prod_{h \in S_{Q}} h\right)^{d \frac{|P|-1}{d}}=(-1)^{\frac{|P|-1}{d} \frac{|Q|-1}{d}}\left(\prod_{\substack{h \in A\\0<|h|< |Q|}}h\right)^{\frac{|P|-1}{d}}.$$
 Substituting into (\ref{(3)}), we get 
  \begin{equation}\label{(4)}
\pi = \left(\left( \prod_{f \in S_{P}} f \right)^{\frac{|Q|-1}{d}}, (-1)^{\frac{|P|-1}{d} \frac{|Q|-1}{d}}\left(\prod_{\substack{h \in A\\0<|h|< |Q|}}h\right)^{\frac{|P|-1}{d}}\right)U.
\end{equation}    
  
Denote $$\widetilde{S_{2}}:=\{h\in A^{+}: 0<|h| < |PQ|~\textrm{and}~(h,PQ)=1\}$$
and $$S_{2}:=\widetilde{S_{2}}\cup g\widetilde{S_{2}}\cup\cdots \cup g^{\frac{q-1}{d}-1}\widetilde{S_{2}}.$$ 
By Chinese reminder theorem, we have  $G_{1}\cong\left(A/(P(t)Q(t))\right)^{*}$ 
and $S_{2}$ be another representation of $U$ in $G$. Since
\begin{equation}
\begin{aligned}
\widetilde{S_{2}} 
    &= \bigl\{f \in A^{+} : 0 < |f| < |P| \bigr\} \\
    &\quad \cup \bigl\{hP + f : h \in A^{+},\; 0 < |h| < |Q|; \; f \in A,\; 0 < |f| < |P| \bigr\} \\
    &\quad \setminus \bigl\{fQ : f \in A^{+},\; 0 < |f| < |P| \bigr\},
\end{aligned}
\end{equation}
 the product of all the elements in $S_{2}$ modulo $P(t)$ is
\begin{equation}
\begin{aligned}
&\quad \frac{\prod_{i=0}^{\frac{q-1}{d}-1}g^{i}\left(\prod_{\substack{h\in A^{+}\\0<|h|<|Q|}}\prod_{\substack{f\in A\\0<|f|<|P|}}(hP+f)\cdot\prod_{\substack{f\in A^{+}\\0< |f|<|P|}}f\right)}
{\prod_{i=0}^{\frac{q-1}{d}-1}g^{i}\prod_{\substack{f\in A^{+}\\0 < |f|<|P|}}(fQ)}\\
&\equiv \frac{\prod_{i=0}^{\frac{q-1}{d}-1}g^{i}\left(\prod_{\substack{f\in A\\ 0< |f|<|P|}}f\right)^{\frac{|Q|-1}{q-1}}}{\prod_{i=0}^{\frac{q-1}{d}-1}g^{i}Q^{\frac{|P|-1}{q-1}}}\\
&\equiv \frac{\left(\prod_{f\in S_{P}}f\right)^{\frac{|Q|-1}{d}}}{\left(\frac{Q}{P}\right)_{d}},
\end{aligned}
\end{equation}
due to (\ref{def}) and the definition of $S_{P}$ (see (\ref{(Sp)})).
Symmetrically, the product modulo $Q(t)$ is $$\frac{\left(\prod_{\substack{h \in A\\0<|h|< |Q|}}h\right)^{\frac{|P|-1}{d}}}{\left(\frac{P}{Q}\right)_{d}}.$$
Thus the product of all the elements in $G$ under the representation $S_{2}$ is
 \begin{equation}\label{(5)}
\pi = \left(\frac{\left(\prod_{f\in S_{P}}f\right)^{\frac{|Q|-1}{d}}}{\left(\frac{Q}{P}\right)_{d}},\frac{\left(\prod_{\substack{h \in A\\0<|h|< |Q|}}h\right)^{\frac{|P|-1}{d}}}{\left(\frac{P}{Q}\right)_{d}}\right)U.
\end{equation} 

Since $U=\langle (\eta,\eta)\rangle=\{(\eta,\eta), (\eta^{2}, \eta^{2}), \cdots, (\eta^{d},\eta^{d})=(1,1)\rangle,$ by comaring (\ref{(4)}) and (\ref{(5)}), we get
\begin{equation}
\begin{aligned}
\left(\frac{P}{Q}\right)_{d}\left(\frac{Q}{P}\right)_{d}^{-1}&=(-1)^{\frac{|P|-1}{d} \frac{|Q|-1}{d}}\\
&=(-1)^{\frac{q-1}{d}\textrm{deg}P\textrm{deg}Q},
\end{aligned}
\end{equation}
which is the desired result.

\end{document}